\documentclass[a4paper,11pt]{article}

\usepackage[utf8]{inputenc}
\usepackage[T1]{fontenc}

\usepackage{amsthm,amsmath}
\usepackage{tikz}
\usepackage{mathrsfs,amssymb,amsfonts} 
\usepackage{enumitem}
\usepackage{fullpage}
\usepackage{hyperref, enumerate}
\usepackage{graphicx}
\usepackage[babel]{microtype}
\usepackage[english]{babel}
\usepackage[capitalise]{cleveref}
\usetikzlibrary{arrows.meta,positioning,shapes.geometric,fit}
\usepackage{thmtools}
\usepackage{mathtools, comment}
\usepackage{amssymb}
\usepackage[nomath]{lmodern}
\usepackage{graphicx}
\usepackage{pgf,tikz,tkz-graph,subcaption}
\usetikzlibrary{arrows,shapes}
\usetikzlibrary{decorations.pathreplacing}
\usepackage{tkz-berge}
\usepackage{enumitem}
\usepackage[normalem]{ulem}
\usepackage{hyperref}
\usepackage{algorithm}
\usepackage{algorithmic}
\hypersetup{colorlinks = true, linkcolor = blue, citecolor = blue, urlcolor = blue}

\usepackage[dvipsnames]{xcolor}
\newcommand{\jj}[1]{{\color{orange} Jorik: #1}}

\newcommand*{\ceilfrac}[2]{\mathopen{}\left\lceil\frac{#1}{#2}\right\rceil\mathclose{}}
\newcommand*{\floorfrac}[2]{\mathopen{}\left\lfloor\frac{#1}{#2}\right\rfloor\mathclose{}}

\allowdisplaybreaks

\usepackage[margin=1in]{geometry}
\newtheorem{defi}{Definition}
\newtheorem{conj}[defi]{Conjecture}

\newtheorem{thr}[defi]{Theorem}
\newtheorem{lem}[defi]{Lemma}

\newtheorem{prop}[defi]{Proposition}

\newtheorem{remark}[defi]{Remark}

\newcommand*{\myproofname}{Proof}

\def\a{\alpha}

\allowdisplaybreaks

\usepackage[margin=1in]{geometry}
\newcommand{\diam}{diam}

\newcommand{\ecc}{ecc}

\title{Maximum total distance under radius and diameter constraints: algorithms, asymptotics, and structure}

\date{}

\author{
Stijn Cambie
\thanks{Department of Computer Science, KU Leuven Campus Kulak-Kortrijk, 8500 Kortrijk, Belgium.
Supported by FWO grant 1225224N.
E-mail: \protect\href{mailto:stijn.cambie@hotmail.com}{\protect\nolinkurl{stijn.cambie@hotmail.com}}}
\and
Jorik Jooken
\thanks{(1) Mathematical Institute, Leiden University, 2333 CA Leiden, the Netherlands, and
(2) Department of Computer Science, KU Leuven Campus Kulak-Kortrijk, 8500 Kortrijk, Belgium.
Supported by FWO grant 1222524N.
E-mail: \protect\href{mailto:jorikjooken@gmail.com}{\protect\nolinkurl{jorikjooken@gmail.com}}}
\and
Mohammad Javad Nadjafi-Arani
\thanks{Faculty of Science, Mahallat Institute of Higher Education, Mahallat, Iran.
E-mail: \protect\href{mailto:mjnajafiarani@gmail.com}{\protect\nolinkurl{mjnajafiarani@gmail.com}} and
\protect\href{mailto:mjnajafiarani@mahallat.ac.ir}{\protect\nolinkurl{mjnajafiarani@mahallat.ac.ir}}}
}

\RequirePackage[normalem]{ulem} 
\RequirePackage{color}\definecolor{RED}{rgb}{1,0,0}\definecolor{BLUE}{rgb}{0,0,1} 
\providecommand{\DIFaddbegin}{} 
\providecommand{\DIFaddend}{} 
\providecommand{\DIFdelbegin}{} 
\providecommand{\DIFdelend}{} 
\providecommand{\DIFaddbeginFL}{} 
\providecommand{\DIFaddendFL}{} 
\providecommand{\DIFdelbeginFL}{} 
\providecommand{\DIFdelendFL}{} 
\newcommand{\DIFscaledelfig}{0.5}
\RequirePackage{settobox} 
\RequirePackage{letltxmacro} 
\newsavebox{\DIFdelgraphicsbox} 
\newlength{\DIFdelgraphicswidth} 
\newlength{\DIFdelgraphicsheight} 
\LetLtxMacro{\DIFOincludegraphics}{\includegraphics} 
\newcommand{\DIFaddincludegraphics}[2][]{{\color{blue}\fbox{\DIFOincludegraphics[#1]{#2}}}} 
\newcommand{\DIFdelincludegraphics}[2][]{
\sbox{\DIFdelgraphicsbox}{\DIFOincludegraphics[#1]{#2}}
\settoboxwidth{\DIFdelgraphicswidth}{\DIFdelgraphicsbox} 
\settoboxtotalheight{\DIFdelgraphicsheight}{\DIFdelgraphicsbox} 
\scalebox{\DIFscaledelfig}{
\parbox[b]{\DIFdelgraphicswidth}{\usebox{\DIFdelgraphicsbox}\\[-\baselineskip] \rule{\DIFdelgraphicswidth}{0em}}\llap{\resizebox{\DIFdelgraphicswidth}{\DIFdelgraphicsheight}{
\setlength{\unitlength}{\DIFdelgraphicswidth}
\begin{picture}(1,1)
\thicklines\linethickness{2pt} 
{\color[rgb]{1,0,0}\put(0,0){\framebox(1,1){}}}
{\color[rgb]{1,0,0}\put(0,0){\line( 1,1){1}}}
{\color[rgb]{1,0,0}\put(0,1){\line(1,-1){1}}}
\end{picture}
}\hspace*{3pt}}} 
} 
\LetLtxMacro{\DIFOaddbegin}{\DIFaddbegin} 
\LetLtxMacro{\DIFOaddend}{\DIFaddend} 
\LetLtxMacro{\DIFOdelbegin}{\DIFdelbegin} 
\LetLtxMacro{\DIFOdelend}{\DIFdelend} 
\DeclareRobustCommand{\DIFaddbegin}{\DIFOaddbegin \let\includegraphics\DIFaddincludegraphics} 
\DeclareRobustCommand{\DIFaddend}{\DIFOaddend \let\includegraphics\DIFOincludegraphics} 
\DeclareRobustCommand{\DIFdelbegin}{\DIFOdelbegin \let\includegraphics\DIFdelincludegraphics} 
\DeclareRobustCommand{\DIFdelend}{\DIFOaddend \let\includegraphics\DIFOincludegraphics} 
\LetLtxMacro{\DIFOaddbeginFL}{\DIFaddbeginFL} 
\LetLtxMacro{\DIFOaddendFL}{\DIFaddendFL} 
\LetLtxMacro{\DIFOdelbeginFL}{\DIFdelbeginFL} 
\LetLtxMacro{\DIFOdelendFL}{\DIFdelendFL} 
\DeclareRobustCommand{\DIFaddbeginFL}{\DIFOaddbeginFL \let\includegraphics\DIFaddincludegraphics} 
\DeclareRobustCommand{\DIFaddendFL}{\DIFOaddendFL \let\includegraphics\DIFOincludegraphics} 
\DeclareRobustCommand{\DIFdelbeginFL}{\DIFOdelbeginFL \let\includegraphics\DIFdelincludegraphics} 
\DeclareRobustCommand{\DIFdelendFL}{\DIFOaddendFL \let\includegraphics\DIFOincludegraphics} 

\begin{document}
\maketitle
\begin{abstract}
    We present an algorithm that determines a graph with the maximum total distance given order $n$ and radius $r$, with time complexity $O(n^2r)$. The same approach immediately yields a modified algorithm for the corresponding problem on trees of given order and diameter, and also gives an algorithm for clique-rooted trees of given order and diameter. For fixed radius or diameter, we also obtain precise asymptotic estimates for the maximum total distance. We investigate the structure of extremal trees and give examples in which branches have different skeletons and unequal orders. We further conjecture that, for fixed diameter $d$ and sufficiently large order $n$, the maximum total distance is attained by a tree when $d$ is even and by a clique-rooted tree when $d$ is odd. Together with our algorithms, this conjecture would yield a practical resolution of the problem of Plesn\'{\i}k~\cite{P84} for sufficiently large order. We obtain such a resolution for the radius variant and for the diameter problem at small order.

\end{abstract}

\section{Introduction}

The \textit{total distance} (also known as the \textit{Wiener index}) of a connected graph $G$ is $W(G)=\sum_{\{u,v\} \in \binom{V(G)}2} d_G(u,v)$, where $d_G(u,v)$ denotes the distance between $u$ and $v$ in $G$. Since the average distance is obtained from $W(G)$ by dividing by $\binom{|V(G)|}{2}$, extremal problems for average distance and total distance are equivalent when the order of the graph is fixed.

One of the most fundamental extremal questions related to the average distance is the (over $40$ years old) question by Plesn\'{\i}k~\cite{P84} to determine the maximum average distance among all graphs, given their order and diameter. For trees, this question also reappeared in other contexts where it is directly related with other quantities, e.g. Kemeny's constant, see~\cite{GW17,JKS23}.
Plesn\'{\i}k's question received a lot of attention and is still not solved exactly in general, although Cambie~\cite{SC19} showed that, for each fixed $d\ge3$, the maximum average distance among graphs of order $n$ and diameter $d$ is $d-\Theta(d^{3/2}/n^{1/2})$ as $n\to\infty$.


The extremal questions with diameter replaced by radius have been investigated as well, both the minimization problem~\cite{CBX13,SC20} and maximization problem~\cite{das2017maximum,SC19}.
In this paper, we study algorithmic versions of the maximization problem. One of our main results is a recursive algorithm that determines the graph with maximum total distance given order $n$ and radius $r$, with time complexity $O(n^2r)$. Every such extremal graph is in fact a tree. The same approach also works for several graph classes where the diameter is given instead of the radius. For trees of given order and diameter, the center is either a vertex or an edge, depending on the parity of the diameter. Therefore, when the diameter is given we can rely on the algorithm that was developed for the case where the radius is given. We also treat the class of \textit{clique-rooted trees} (a clique-rooted tree is a clique with trees attached to its vertices). This class is particularly relevant for the case where the diameter is odd, since we expect that in such cases there is a central clique rather than a central vertex for an extremal graph. We remark that previous computational approaches were only able to determine the extremal trees for a small number of vertices (say, 20), whereas our algorithms can easily determine the extremal trees having thousands of vertices. A related problem on maximizing the Wiener index of a tree, given the degree sequence, was studied with an algoritmic point of view in~\cite{CSWW11}.

Our results are motivated by the conjectures by Sun et al.~\cite[Conj.~3 \& 5]{SIRV19} and~\cite[Ques.~6.1]{SC19}.
 The available evidence suggests the following generalization of their conjectures. 

\begin{conj}\label{conj:Nlarge}
    For $n$ sufficiently large in terms of $d$, the maximum total distance among all graphs of order $n$ and diameter $d$ is obtained by 
    \begin{itemize}
        \item a tree if $d$ is even,
        \item a clique-rooted tree if $d$ is odd.
    \end{itemize}
\end{conj}

A further question is how large $n$ must be, in terms of $d$,
for this conjectured behaviour to hold.
For $d=5$, it turns out that for $n=12$ already, the extremal graph is a tree (with $W=198$), with the cycle $C_{11}$ with one vertex duplicated as second best edge-minimal example ($W=197$).  For $d=5$, we conjecture that an extremal graph can be chosen to be a tree for $6\le n\le20$, except when $n=11$, and can be chosen to be a clique-rooted tree with a clique of order at least $3$
for $n>20$. The value $n=2d+1$ deserves a special mention due to the following conjecture by DeLaVi{\~n}a and Waller~\cite{DLVW08}, proven in~\cite{Cam26+}.

\begin{conj}[\cite{DLVW08}]\label{conj:DLVW}
    Let $G$ be a graph with diameter $d>2$ and order $2d+1$. Then $W(G) \leq W(C_{2d+1})$.
\end{conj}

For large values of $d$, the regime of exceptions around $2d+1$ is expected to be larger.

The rest of this paper is organized as follows. In \cref{subsec:not&def}, we introduce the terminology and notation used throughout the paper. Next, in \cref{sec:prelim}, we state preliminary results. In~\cref{sec:algo}, the algorithms are presented to compute the extremal trees and clique-rooted trees. In particular, we completely resolve the algorithmic version of~\cite[Prob.~2]{CBX13}, which asks the maximum total distance given order and radius, also popularised in e.g.~\cite[Prob.~5.3]{KST16} and~\cite[Prob.~5.9]{KST24}. 
For fixed radius $r$, in~\cref{sec:given_r}, we prove that $W(G) \le rn^2- 2\sqrt{\mu_r}n^{3/2}+O_r(n)$, where $\mu_r$ can be efficiently computed exactly for moderate $r$ (say, up to $100$) and obtain the approximate form of the branches for large order. In~\cref{sec:structure}, we discuss the structure of extremal graphs of given order and diameter, and provide examples of extremal trees that are quite unbalanced. 
The results convince us that (in general) a clean characterization of the extremal trees of given order and radius does not exist and therefore algorithms are necessary to find such trees, while we understand the rough structure.
Finally, in~\cref{sec:conc} we recap with the expected answer for the problem by Plesn\'{\i}k.

\subsection{Terminology and notation}\label{subsec:not&def}

The open and closed neighbourhood of a subset $U \subset V$ are denoted by $N(U)$ and $N[U]$. 
In a tree \(T = (V, E)\), for any set of vertices \(X \subseteq V\), let \(X_i = N_i(X)\) be the set of vertices at distance exactly \(i\) from \(X\). In a special case, if $L$ is the set of leaves, $L_i$ denotes the set of vertices for which the closest leaf is at distance (exactly) $i$. In particular, $L_1=N(L).$

The eccentricity of a vertex $v$, $\ecc(v)$, equals $\max_{u \in V} d(u,v).$
The radius and diameter are, respectively, $\min_v \ecc (v)$ and $\max_v \ecc (v)$. A \textit{clique-rooted tree} is a graph obtained from a forest by adding exactly all edges between the roots of the trees so that the roots of all trees induce a clique. We will also use vertex-weighted graphs for our arguments. For a vertex-weighted graph $G$, where every vertex $v$ has a weight (multiplicity) $w(v)$, the total distance can be defined as $$W(G)= \sum_{\{u,v\} \in \binom{V(G)}{2} } w(u)w(v)d(u,v).$$
Note that the unweighted case corresponds to taking $w(v)=1$ for all vertices $v$.

\section{Preliminary results}
\label{sec:prelim}
We start with mentioning two results from the literature. 

\begin{lem}\label{lem:nr_tree}[\cite[Lem.~3]{CBX13}]
    Any graph maximizing the total distance among all graphs with order $n$ and radius $r$ is a tree.
\end{lem}

\begin{lem}\label{lem:leaves_far}[\cite[Lem.~4.1]{SC19}, \cite[Lem.~2.1]{BDGSSZ25}]
A tree with maximum total distance among all trees of order $n$ and diameter $d$ does not contain any leaf that has distance to a central vertex strictly smaller than $\frac{d-1}{2}$, equivalently $\ecc{(\ell)}=d$ for every leaf $\ell.$
\end{lem}

The same statement is true for clique-rooted trees with given order and diameter. In an extremal clique-rooted tree with a clique of order at least three, the maximum height $H$ of the attached rooted trees must be attained at least twice. Indeed, otherwise, denoting the second-largest height by $h$, replacing the clique by a star centred at the root of the tallest tree strictly increases the total distance without increasing the diameter, since $2h+2\leq H+h+1\leq d$. Hence $d=2H+1$. If a rooted tree attached to the clique has height less than $H$, we can delete all but one of the clique edges incident with its root. This strictly increases the total distance while preserving the diameter, so all rooted trees attached to the clique have height $H$.

Suppose now that a leaf $\ell$ is at depth less than $H$, and let $v$ be its closest ancestor with at least two children. Let $P$ be the pendant path from $v$ to $\ell$, excluding $v$, let $w$ be a child of $v$ outside $P$, and let $u$ be the parent of $v$ (or another clique vertex if $v$ belongs to the clique). Reattaching $P$ to either $u$ or $w$ preserves the diameter, and at least one of these operations strictly increases the total distance, since for every vertex $x$ outside $P$ we have $d(u,x)+d(w,x)\geq 2d(v,x)$,
with strict inequality for $x=v$.

The following monotonicity lemma is useful when going from the setting of given radius to given diameter.

\begin{lem}
	\label{lem:monotonicity}
		Let $W^T_{n,d}$ be the maximum total distance of any tree with diameter $d$ and order $n$.
		Then for any $d\ge 2$ and $n\ge d+2$, we have $W^T_{n,d+1}>W^T_{n,d}.$ 
\end{lem}

\begin{proof}
    Let $T$ be an extremal tree of order $n$ and diameter $d$.
    By~\cref{lem:leaves_far}, all leaves are equally far from the centers.
    Now let $v$ be a vertex of $T$ (which is not a path) with degree at least $3$ which is nearest to a leaf among all such vertices.
    Now by this choice, $v$ is equally far (say distance $k$) from (at least) two leaves $\ell$ and $m$.
    Now replace the paths between $v$ and $\ell$ or between $v$ and $m$ by paths with length $k-1$ and $k+1$ respectively (equivalently, delete the edge incident to $m$ and add the edge $\ell m$).
    Then it is immediate to see that for every other vertex, the sum of distances towards the $2k+1$ vertices on the union of the two paths increased and so $W$ increased while the diameter became $d+1.$
\end{proof}

Next, we prove a lemma that most branching vertices have degree at most $3$.
	
		\begin{lem}
		\label{lem:deg2or3}
		Let $T=(V,E)$ be a tree attaining the maximum total distance $W_{n,d}^T$ among all trees with diameter $d$ and order $n$.
        Let $L_0=L$ be the set of leaves,  and let $C$ be the center (or set of centers) of $T$.
		The degree of any vertex $u \in V \setminus (N[L] \cup C)$ is bounded by $3.$
		\end{lem}

	\begin{proof}
Assume, for contradiction, that the claim is false. 
Let $v$ be a center of $T$ that is closest to $u$, and suppose $u \in N_i(v)$ for some $i \ge 1$. 
Let $p \in N_{i-1}(v)$ be the neighbour of $u$ on the unique path from $u$ to $v$.

The vertex $u$ has at least three neighbours in $N_{i+1}(v)$, each of which serves as the root of a distinct subtree. 
Among these, let $T_2$ and $T_3$ denote the two smallest such subtrees (with respect to the number of vertices), and let $T_1$ be the union of all remaining subtrees attached to $u$ via neighbours in $N_{i+1}(v)$.
A schematic illustration is provided in \cref{fig:1}. 
Note that in general, $T_1$ may consist of multiple neighbours of $u$ (i.e., more than one subtree rooted at distinct vertices in $N_{i+1}(v)$).

\begin{figure}
    \begin{center}
    	\begin{tikzpicture}[
		>=Stealth,
		every node/.style={font=\large},
		vertex/.style={circle,draw,minimum size=6mm,inner sep=1pt},
		cloud/.style={draw,ellipse,minimum width=2.2cm,minimum height=1.1cm,align=center}
		]
		
		\node (v) [vertex] at (-5,0){$v$};
		\node (dots) at (-4,0) {$\cdots$};
		\node (p)   [vertex] at (-2.8,0) {$p$};
		\node (u)   [vertex] at (-1,0) {$u$};
		
		\draw (p) -- (u);
		
		\node (v1) [vertex] at (.8,1.4) {$v_1$};
		\node (v2) [vertex] at (.8,0) {$v_2$};
		\node (v3) [vertex] at (.8,-1.4) {$v_3$};
		
		\draw (u) -- (v1);
		\draw (u) -- (v2);
		\draw (u) -- (v3);
		
		\draw (2.2,1.4) ellipse (1.1cm and 0.55cm);
		\node at (2.3,1.4) {$T_1$};
		
		\draw (2.2,0) ellipse (1.1cm and 0.55cm);
		\node at (2.3,0) {$T_2$};
		
		\draw (2.2,-1.4) ellipse (1.1cm and 0.55cm);
		\node at (2.3,-1.5) {$T_3$};
			
	\end{tikzpicture}
    		\caption{Part of  $T$} \label{fig:1}
\end{center}
\end{figure}
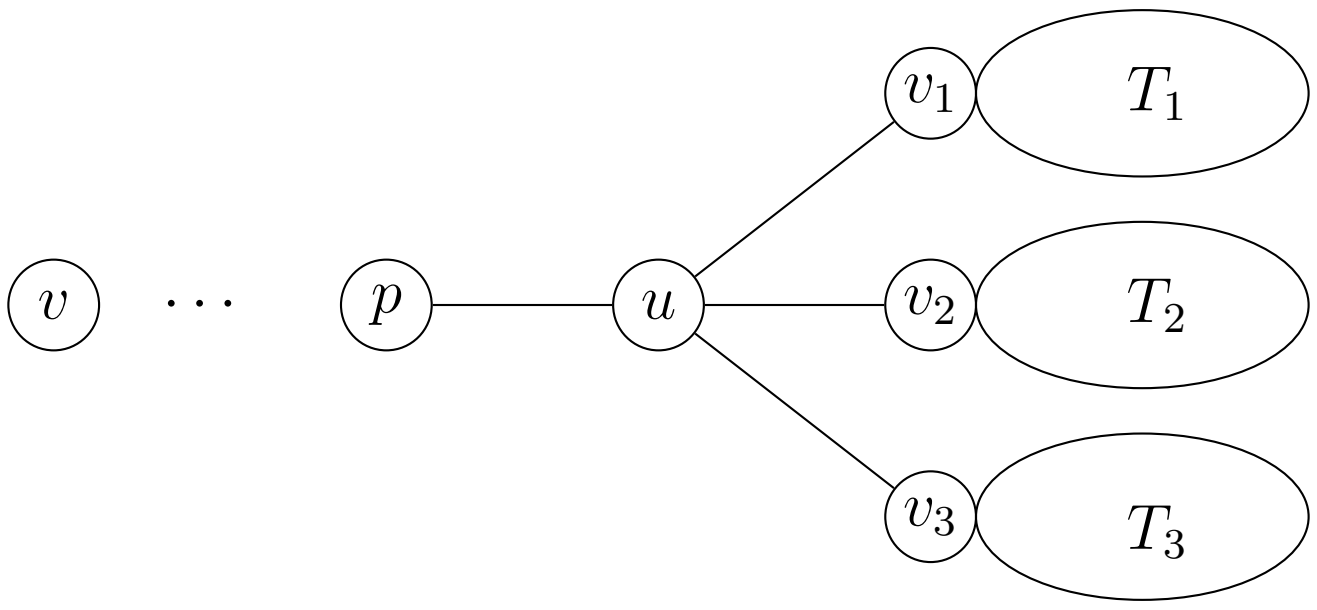









	We identify each $T_j$ with its vertex set, and let $v_2$ and $v_3$ be the roots of $T_2$ and $T_3$, respectively. Let $X_1=N_T(u)\cap T_1$ be the set of roots of the subtrees comprising $T_1$. In both constructions, first move all subtrees rooted at the children of $v_3$ to $v_2$: for each child $x$ of $v_3$, delete $v_3x$ and add $v_2x$. Next, to obtain $T'$, delete $uv_3$ and add $v_2v_3$ (\Cref{fig:2}). To obtain $T''$, instead delete $uv_3$ and add $pv_3$, and for every $x\in X_1$, delete $ux$ and add $v_3x$. These operations are depicted in \Cref{fig:3}. 
	
	Both modified graphs are trees. Root the original tree towards its center (or towards the nearer endpoint of its central edge). The descendants originally below $v_3$ retain their depths in both constructions, and all vertices of $T_1$ retain their depths in $T''$. The only vertex whose depth increases is $v_3$ in $T'$, and its new depth is at most the original maximum depth, since $v_3$ was not a leaf. Thus neither construction increases the maximum depth on either side of the center. By Lemma~\ref{lem:leaves_far} and the choice of $u$, a leaf at the original maximum depth in the modified part remains at that depth, for instance a leaf in $T_2$. A leaf at maximum depth in another central branch (or on the opposite side of the central edge) is unchanged. Their distance is still $d$. Consequently, $\diam(T')=\diam(T'')=d$.
		
	The diameter of both $T'$ and $T''$ is equal to the diameter of $T$.
	

    \begin{figure}[ht]
			\begin{minipage}[b]{.49\linewidth}
			\begin{center}
		
		\begin{tikzpicture}[
		>=Stealth,
		every node/.style={font=\large},
		vertex/.style={circle,draw,minimum size=5mm,inner sep=1pt},
		cloud/.style={draw,ellipse,minimum width=2cm,minimum height=1.1cm,align=center}
		]
		
		\node (p)   [vertex] at (-2.8,0) {$p$};
		\node (u)   [vertex] at (-1,0) {$u$};
		
		\draw (p) -- (u);
		
		\node (v1) [vertex] at (.8,1.4) {$v_1$};
		\node (v2) [vertex] at (.8,0) {$v_2$};
		\node (v3) [vertex] at (2.3,-1.5) {$v_3$};
		
		\draw (u) -- (v1);
		\draw (u) -- (v2);
		\draw (v2) -- (v3);
		
		\draw (2.15,1.4) ellipse (1.1cm and 0.55cm);
		\node at (2.3,1.4) {$T_1$};
		
		\draw (2.15,0) ellipse (1.1cm and 0.55cm);
		\node at (2.3,0) {$T_2$};
		
		\draw (.8,-1.35) ellipse (.55cm and 1.1cm);
		\node at (.8,-1.5) {$T_3$};
		
		
	\end{tikzpicture}
	\end{center}
	\subcaption{part of  $T'$} \label{fig:2}
\end{minipage}
\begin{minipage}[b]{.49\linewidth}
\begin{center}
			\begin{tikzpicture}[
		>=Stealth,
		every node/.style={font=\large},
		vertex/.style={circle,draw,minimum size=5mm,inner sep=1pt},
		cloud/.style={draw,ellipse,minimum width=2cm,minimum height=1.1cm,align=center}
		]
		
		\node (p)   [vertex] at (-2.8,0) {$p$};
		\node (u)   [vertex] at (-1,0) {$u$};
		
		\draw (p) -- (u);
		
		\node (v1) [vertex] at (.8,1.4) {$v_1$};
		\node (v2) [vertex] at (.8,0) {$v_2$};
		\node (v3) [vertex] at (-1,1.4) {$v_3$};
		
		\draw (v3) -- (v1);
		\draw (u) -- (v2);
		\draw (p) -- (v3);
		
		\draw (2.15,1.4) ellipse (1.1cm and 0.55cm);
		\node at (2.3,1.4) {$T_1$};
		
		\draw (2.15,0) ellipse (1.1cm and 0.55cm);
		\node at (2.3,0) {$T_2$};
		
		\draw (.8,-1.35) ellipse (.55cm and 1.1cm);
		\node at (.8,-1.5) {$T_3$};
		
		
	\end{tikzpicture}
\end{center}
\subcaption{part of $T''$} \label{fig:3}
\end{minipage}
\caption{Partial illustration of two modified versions of $T$}
\end{figure}
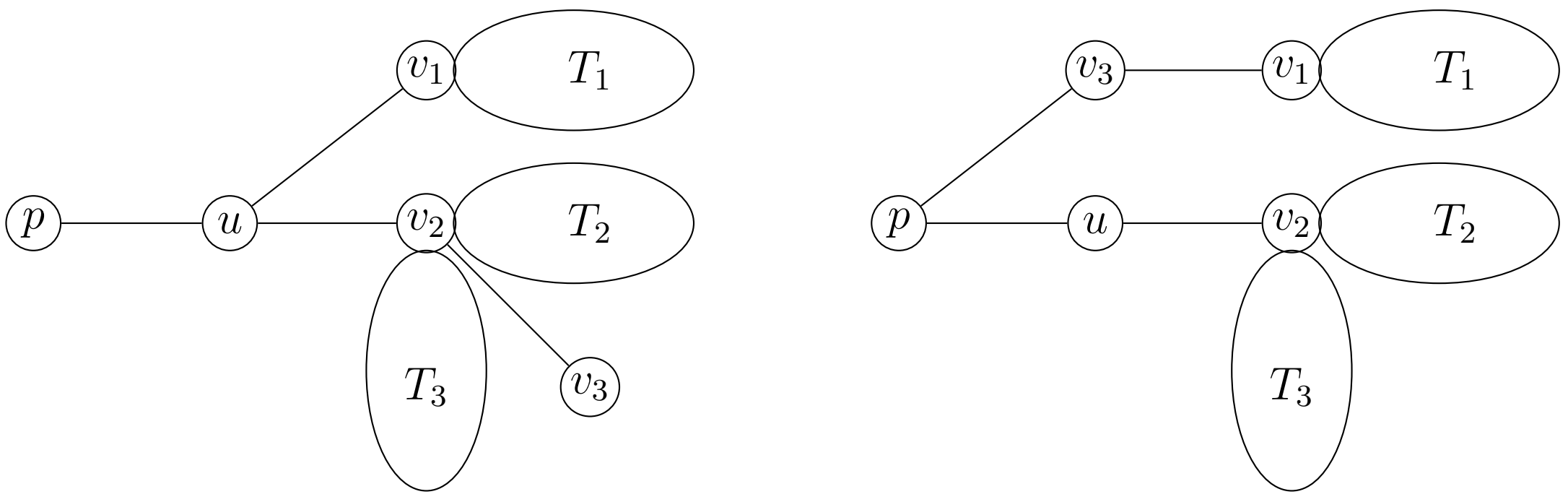
Now the pairs of vertices $(w_1,w_2)$ for which $d_T(w_1,w_2) \neq d_{T'}(w_1,w_2)$ are the following ones:
	
	$$\begin{array}{ll}
	d_T(w_1,w_2) = d_{T'}(w_1,w_2)+2 & w_1 \in T_2, w_2 \in T_3 \backslash \{v_3\},\\
	d_T(v_3,w) = d_{T'}(v_3,w) - 1 & w\in V(T)\setminus(T_2\cup\{v_3\}),\\
	d_T(v_3,w) = d_{T'}(v_3,w) + 1 & w \in T_2,
	\end{array}$$
	and hence, 
	$$W(T) - W(T') =  2|T_2||T_3|+1-n.$$
	
Similarly, the pairs of vertices $(w_1,w_2)$ for which $d_T(w_1,w_2) \neq d_{T''}(w_1,w_2)$ are the following ones:
	
	$$\begin{array}{ll}
	d_T(w_1,w_2)=d_{T''}(w_1,w_2) - 2 & w_1 \in T_1, w_2 \in T_2 \cup (T_3 \backslash \{v_3\}) \cup \{u\},\\
	d_T(w_1,w_2)=d_{T''}(w_1,w_2) + 2 & w_1 \in T_2, w_2 \in T_3\backslash \{v_3\},\\
	d_T(v_3,w)=d_{T''}(v_3,w) + 1  & w\in V(T)\setminus(T_2\cup T_3\cup\{u\}),\\
	d_T(v_3,w)=d_{T''}(v_3,w) - 1 &  w \in T_2 \cup \{u\},\\
	d_T(v_3,w)=d_{T''}(v_3,w) - 3 & w \in T_3\backslash \{v_3\},  
	\end{array}$$
	and hence,
	\begin{eqnarray*}
		W(T) - W(T'') &=& -2|T_1|(|T_2| + |T_3|) + 2|T_2|(|T_3|-1) \\ &+& (n - |T_2| - |T_3|-1) - (|T_2| + 1) - 3(|T_3| - 1) \\ &=&
		n + 1 + 2|T_2||T_3| - 2|T_1||T_2| - 2|T_1||T_3| - 4|T_2| - 4|T_3|.
	\end{eqnarray*}  
	
	Since, $T$ attains the maximum Wiener index, we need $W(T) - W(T') \geq 0$ and $W(T) - W(T'') \geq 0$,
	thus the sum of two has to satisfy
	$$ 4|T_2||T_3| - 2|T_1||T_2| - 2|T_1||T_3| - 4|T_2| - 4|T_3| + 2 \geq 0.$$
	But this is not the case since $\lvert T_1 \rvert \ge |T_2|,|T_3|>1$.
\end{proof}

The statement of \cref{lem:deg2or3} also holds for extremal clique-rooted trees, with $C$ denoting the vertex set of the central clique, since both transformations still preserve the diameter and the corresponding distance calculations remain valid.

\section{Algorithm(s) computing maximum total distance given n and $d$ or $r$}\label{sec:algo}

\subsection{The recursive functions}\label{subsec:functions}

For trees, we can iteratively compute the total distance by extending our computations to weighted trees, and noticing that $W(T)$ can also be computed by summing the number of times every edge is traversed (technically; the Wiener index and Szeged index are equal for (weighted) trees).
If one denotes for an edge $e=uv$ with $n_u(e)$ and $n_v(e)$ the total weight of the components (containing $u$ and $v$ resp.) of the tree after deleting $e,$ we have $W(T)= \sum_{e=uv} n_u(e)n_v(e)$.
This alternative computation motivates and explains the following help functions, that are also useful in general (weighted) graphs which contain subtrees.

\begin{defi}\label{def7}
     Fix the total order $n$. Let $dW(n,t,i)$ be the maximum weighted Wiener index of a rooted tree of height exactly $i$ with $t$ non-root vertices, each of weight one, and a root of weight $n-t$. All leaves (in the rooted sense) are at height $i$, and every vertex at depth strictly less than $i-1$ has at most two children. Vertices at depth $i-1$ may have arbitrarily many children. For $i=0$, the tree consists only of the root. We set
	\[
	dW(n,0,0)=0,\qquad dW(n,t,0)=-\infty\quad(t>0),
	\]
	and give every infeasible state the value $-\infty$. Throughout, $0\leq t\leq n-1$. Let $dW_1(n,t,i)$ denote the corresponding maximum under the additional condition that the root has exactly one child, for $i\geq1$.

\end{defi}

        \begin{prop}\label{prop:dW_iteration}
             For $1\leq t\leq n-1$,
	\[
	dW(n,t,1)=t(n-1).
	\]
	For $i\geq1$ and $1\leq t\leq n-1$,
	\[
	dW_1(n,t,i)=dW(n,t-1,i-1)+t(n-t).
	\]
	For $i\geq2$,
	\begin{align*}
		dW_2(n,t,i)&=\max_{\substack{a+b=t\\a,b\geq1}}
		\bigl\{dW_1(n,a,i)+dW_1(n,b,i)\bigr\},\\
		dW(n,t,i)&=\max\{dW_1(n,t,i),dW_2(n,t,i)\}.
	\end{align*}
	Here $dW_2$ is the maximum restricted-state value with two children at the root, and a maximum over an empty set is $-\infty$. In particular, $dW_1(n,1,1)=n-1$ and $dW_1(n,t,1)=-\infty$ for $t\ne1$. For positive heights, states with $t=0$ are infeasible.
	        \end{prop}

	\begin{proof}
		When $i=1$, the only possible rooted tree is a star with $t$ leaves, and each pendant edge contributes $n-1$ to the weighted Wiener index. If the root has one child, the edge incident with the root separates $t$ unit-weight vertices from the remaining weight $n-t$, contributing $t(n-t)$. Contracting this edge gives the state $dW(n,t-1,i-1)$, since the new root has weight $n-t+1$. Conversely, an optimal solution to that smaller state can be extended by reinstating the edge.
		
		If the root has two children, the components below them have orders $a$ and $b$, with $a+b=t$. By the edge formula for the Wiener index, their contributions, including their attachment edges, are additive. They are precisely the two single-branch states in the displayed recurrence. Conversely, independently optimal realizations of these states can be attached to the same root, without violating the height or branching restrictions. For $i\geq2$, the root has either one or two children, proving the last equality.
	\end{proof}

	For a central root the number of children need not be bounded by two. We therefore introduce a separate function $C(n,t,h)$, defined as the maximum contribution of the edges in a rooted tree of height $h$ with $t$ non-root vertices, in which every branch at the root satisfies the restrictions of Definition \ref{def7}. For $h\geq1$ and $t\geq1$,
	\[
	C(n,t,h)=\max_{\substack{s\geq1,\ t_j\geq h\\\sum_{j=1}^s t_j=t}}
	\sum_{j=1}^sdW_1(n,t_j,h).
	\]
	For the purpose of the exact-fill recurrence below, the empty collection of branches has value zero. Thus $C(n,0,h)=0$ is an auxiliary initialization, not a rooted tree of positive height. Set $C(n,0,0)=0$ and $C(n,t,0)=-\infty$ for $t>0$.
\begin{defi}
  For a rooted component of order $t$ incident with the central edge (or a central clique vertex in the setting under consideration), define
	\[
	F(n,t,r):=C(n,t-1,r-1)+\frac{t(n-t)}{2}.
	\]
	The factor $1/2$ accounts for the contribution of the central edge or the weighted clique when the individual rooted components are combined.
  
\end{defi}

\subsection{The algorithm}

We now discuss an algorithm to compute the maximum total distance among all
graphs of order $n$ and radius $r$\footnote{We propose an algorithm based on dynamic programming. By keeping track of which choices are optimal, one can also compute the optimal graph itself instead of only the optimal value.}. By \cref{lem:nr_tree}, the maximum is attained by a tree. If $n=2r$, there is only one such tree (namely the path $P_{2r}$), so we assume that $n \geq 2r+1$. By \cref{lem:monotonicity}, the maximum total distance among trees of radius $r$ is attained by a tree of diameter $2r$. Hence, this tree has a unique center $c$.

Let $T$ be this tree and let the components of $T-c$ have orders $t_1,t_2,\ldots,t_s$. We have $\sum_{i=1}^{s}t_i=n-1$. By \cref{lem:leaves_far}, every connected component together with its edge to the center $c$, is a rooted branch of height $r$. By \cref{lem:deg2or3} and \cref{prop:dW_iteration}, the maximum contribution of such a branch of order $t_j$ is $dW_1(n,t_j,r)$. Therefore, the maximum total distance among all graphs of order $n$ and radius $r$ is:
$$
    \max_{\substack{s \geq 2\\
    t_1,\ldots,t_s \geq r\\
    \sum_{j=1}^{s}t_j=n-1}}
    \sum_{j=1}^{s} dW_1(n,t_j,r).
$$

This can in fact be seen as a special instance of the \textit{unbounded knapsack problem}~\cite{kellerer2004unbounded}: given a knapsack having a capacity $c$ and $n$ item types, where each item type contains arbitrarily many identical items with weight $w_i$ and profit $p_i$ ($i \in [n]$), the goal is to select a set of items such that the sum of the weights is bounded by $c$ and the sum of the profits is maximized. In our case, the capacity $c$ is equal to $n-1$ and for each $i$ between $r$ and $n-1-r$, there is an item type where every item has weight $i$ and profit $dW_1(n,i,r)$. Moreover, we require the capacity to be filled exactly. The upper bound $n-1-r$ on the item weights ensures that at least two branches are selected. The exact-fill constraint can be incorporated into the standard dynamic programming algorithm by initializing all positive capacities with profit $-\infty$ instead of zero.

This unbounded knapsack problem can subsequently be translated into a
\textit{$0$-$1$ knapsack problem} (where there is exactly one item available of each item type) by replacing each item type by a logarithmic number of items. More precisely, for an item type of weight $i$, any feasible solution uses at most $\left\lfloor \frac{n-1}{i}\right\rfloor$ items of this type. The idea is to replace this item type by an item with weight $2^h i$ and profit $2^h dW_1(n,i,r)$ for increasingly larger integers $h$ as long as $2^h \leq \left\lfloor\frac{n-1}{i}\right\rfloor$ (starting with $h=0$). By design, every choice of $m \leq \left\lfloor\frac{n-1}{i}\right\rfloor$ copies of the original item type $i$ (in the unbounded knapsack problem) can then be represented by a subset of these new items. Moreover, any subset representing more than $\left\lfloor\frac{n-1}{i}\right\rfloor$ copies has weight larger than $n-1$, yielding an infeasible solution.

We now show that the number of items in the resulting $0$-$1$ knapsack instance is at most $2n$. For a fixed value of $h$, an item with weight $2^h i$ is introduced only if $2^h \leq \left\lfloor\frac{n-1}{i}\right\rfloor$, i.e., only if $i \leq \frac{n-1}{2^h}$. Hence, for this fixed $h$ at most $\left\lfloor\frac{n-1}{2^h}\right\rfloor$ items are introduced. Thus, the total number of items is at most
$$
\sum_{h\geq 0}\left\lfloor\frac{n-1}{2^h}\right\rfloor
    \leq \sum_{h\geq 0}\frac{n-1}{2^h}<2n.
$$
The capacity of the knapsack is equal to $n-1$ and therefore this $0$-$1$ knapsack instance can be solved in $O(n^2)$ time using a standard dynamic programming algorithm~\cite{kellerer2004exact}. The values from~\cref{prop:dW_iteration} can be computed in $O(n^2r)$ time using the given recursions, leading to a total time of $O(n^2r)$ for the entire algorithm.

We next derive a related expression for clique-rooted trees of diameter $2r-1$.
By the observations in Section~2, we may restrict ourselves to a central clique of order at least 2 with rooted trees attached to its vertices, each of height $r-1$ and with all leaves at that height.
If the extremal graph is a tree, we take its central edge as the central clique. Let $T$ be an optimal graph of this form, let $c_1,\ldots,c_s$ be the vertices of the central clique of $T$ and let $t_j$ be the order of the rooted tree attached to $c_j$ (including $c_j$ itself). We have $\sum_{j=1}^{s}t_j=n$.

The maximum contribution to the total distance of the edges inside the rooted tree of order $t_j$ is $C(n,t_j-1,r-1)$. Moreover, the total contribution to the total distance for the edges inside the central clique is
$$
\sum_{1\leq a<b\leq s} t_a t_b
=
\frac 12 \sum_{a=1}^{s} t_a \sum_{b \not =a} t_b
=
\sum_{j=1}^{s}\frac{t_j(n-t_j)}{2}.
$$

Therefore, the maximum total distance among graphs of this form and order $n$ is given by:
$$
    \max_{\substack{s\geq 2\\
    t_1,\ldots,t_s\geq r\\
    \sum_{j=1}^{s}t_j=n}}
    \sum_{j=1}^{s} F(n,t_j,r).
$$

This is again an unbounded knapsack problem, which can be translated into a $0$-$1$ knapsack problem solvable in $O(n^2)$ time.\footnote{The above algorithms were also implemented and made publicly available at\\ \url{https://github.com/JorikJooken/maximizeTotalDistanceGivenOrderAndRadius}.}

Finally, we discuss how the corresponding problems with given diameter $d$ can be handled. If $d=2r$, then every tree of diameter $d$ has radius $r$ and a unique center. Therefore, in this case the first algorithm can be used with radius $r=d/2$.

If $d=2r-1$, then every tree of diameter $d$ has a central edge (i.e., a special case of the clique-rooted setting in which the center is a $K_2$). Hence, it suffices to try all possible orders of one of the two branches. If one branch has order $t$, then the other one has order $n-t$. Therefore, the maximum total distance among trees of order $n$ and diameter $d=2r-1$ is
$$
    \max_{r\leq t\leq n-r}
    \left(
        F(n,t,r)+F(n,n-t,r)
    \right).
$$
This is precisely the clique-rooted computation restricted to two rooted branches, and it can clearly be computed in $O(n)$ time once the values of $F$ have been computed. If the central clique is not restricted to consist of two vertices, the clique-rooted computation above gives the optimum among clique-rooted trees of diameter $d=2r-1$.

Hence, we obtain exact algorithms for the corresponding optimization problems over trees and over clique-rooted trees for given diameter $d$. By \cref{conj:Nlarge}, an extremal graph among all graphs of order $n$ and diameter $d$ can be chosen from the appropriate one of these classes (when $n$ is large enough in terms of $d$). Under this conjecture, the same algorithms therefore also compute the maximum total distance among all graphs of order $n$ and diameter $d$.
\section{Precise asymptotic estimate given $r$}\label{sec:given_r}


Next, we determine the leading correction term in the asymptotic
estimates, as such extending the results by Mukwembi and Vetr\'{\i}k~\cite{MV14}. Throughout this section,
$r$ is fixed and $n$ tends to infinity.

\begin{defi}\label{def:mu_r}
For $r\ge 2$, an \emph{$r$-branch skeleton} is a rooted tree $S$ whose root
is at level $1$, whose terminal vertices are precisely the vertices at level
$r-1$, and in which every nonterminal vertex has one or two children. The
terminal vertices will eventually receive pendant leaves, which then lie at
level $r$. Let $Z(S)$ denote the set of terminal vertices.

Adjoin a new vertex $o$ at level $0$, adjacent to the root of $S$, and denote
the resulting tree by $\widehat S$. For $v\in V(S)$, let $\ell(v)$ be its
level and define
\[
    b(S):=\sum_{v\in V(S)}\bigl(r-\ell(v)\bigr).
\]
Thus $b(S)$ is the total deficit in depth of the vertices of the skeleton
relative to level $r$.

For an edge $e\in E(\widehat S)$, let $Z_e$ be the set of terminals in the
component of $\widehat S-e$ not containing $o$. Define the (effective) resistance\footnote{This is the usual definition known by Thomson’s principle when all vertices of $Z(S)$ are identified}
\[
    R(S):=
    \min_{\substack{\lambda_z\ge 0\ (z\in Z(S))\\
                    \sum_{z\in Z(S)}\lambda_z=1}}
    \sum_{e\in E(\widehat S)}
    \left(\sum_{z\in Z_e}\lambda_z\right)^2.
\]

Finally, put
\[
    \mu_r:=
    \min\bigl\{b(S)R(S):S\text{ is an $r$-branch skeleton}\bigr\}.
\]
\end{defi}

For fixed $r$, there are only finitely many $r$-branch skeletons, and hence
the minimum in Definition~\ref{def:mu_r} is well-defined.

\begin{thr}\label{thr:fixed-diameter-asymptotics}
For every fixed integer $r\ge 2$,
\[
    \max_{\substack{|V(T)|=n\\ \diam(T)=2r}}W(T)
    =
    rn^2-2\sqrt{\mu_r}\,n^{3/2}+O_r(n),
\]
and
\[
    \max_{\substack{|V(T)|=n\\ \diam(T)=2r+1}}W(T)
    =
    \left(r+\frac14\right)n^2
    -2\sqrt{\mu_r}\,n^{3/2}+O_r(n),
\]
where both maxima are taken over trees.
\end{thr}

\begin{proof}
By~\cref{lem:leaves_far},
in a Wiener-maximal tree of the prescribed diameter, all leaves have the same distance towards the center. Moreover, by~\cref{lem:deg2or3}, after the leaves are
removed, each component away from the centre or central edge is an
$r$-branch skeleton. Thus a branch with skeleton $S$ is obtained by
attaching $a_z\ge 1$ leaves to each $z\in Z(S)$.

We first estimate the contribution of a single branch $B$. Let $c$ be the central vertex to which the branch is attached. 
Put
\[
    A:=\sum_{z\in Z(S)}a_z,
    \qquad
    t:=|V(S)|+A.
\]
For an edge $e$, let $s_e$ be the order of the component of $T-e$
not containing $c$. Since $\sum_e s_e$ is the sum of the levels of all
vertices in the branch, we have
\begin{equation}\label{eq:branch-linear-term}
    \sum_e s_e
    =\sum_{v\in V(B)}d_T(c,v)
    =rt-b(S).
\end{equation}

For $e\in E(\widehat S)$, let $q_e$ be the number of vertices of $S$ in the
component below $e$ (the part of $S \setminus e$ not containing the center $o$). Every pendant edge has $s_e=1$, and therefore
\begin{equation}\label{eq:branch-quadratic-term}
    \Phi(B):=\sum_e s_e^2
    =
    A+
    \sum_{e\in E(\widehat S)}
    \left(q_e+\sum_{z\in Z_e}a_z\right)^2.
\end{equation}

Writing $\lambda_z=a_z/A$ and omitting the nonnegative terms involving the
$q_e$, we obtain
\[
    \Phi(B)
    \ge
    A^2
    \sum_{e\in E(\widehat S)}
    \left(\sum_{z\in Z_e}\lambda_z\right)^2
    \ge A^2R(S).
\]
Since $A=t-|V(S)|$ and $|V(S)|=O_r(1)$, this gives
\[
    \Phi(B)\ge R(S)t^2-O_r(t).
\]

Conversely, let $(\lambda_z)_{z\in Z(S)}$ be a minimizing vector in the
definition of $R(S)$. For any sufficiently large $A$, choose positive
integers $a_z$ with sum $A$ such that
\[
    a_z=A\lambda_z+O_r(1)
    \qquad (z\in Z(S)).
\]
Substitution in \eqref{eq:branch-quadratic-term} gives
\[
    \Phi(B)\le R(S)A^2+O_r(A)
             =R(S)t^2+O_r(t).
\]
Consequently, uniformly over all $r$-branch skeletons,
\begin{equation}\label{eq:branch-energy-asymptotic}
    \min_{\substack{B\text{ has skeleton }S\\ |V(B)|=t}}
    \Phi(B)
    =
    R(S)t^2+O_r(t).
\end{equation}

We first consider diameter $2r$. Let $c$ be the centre of $T$, and let
$B_1,\ldots,B_k$ be its branches. Write $t_i$ for the order of $B_i$ and
put
\[
    b_i:=b(S_i),
    \qquad
    R_i:=R(S_i),
\]
where $S_i$ is the skeleton of $B_i$. By the edge formula
\[
    W(T)=\sum_{e\in E(T)}s_e(n-s_e)=n \sum_{e\in E(T)}s_e - \sum_{i=1}^k \Phi(B_i),
\]
together with \eqref{eq:branch-linear-term} and
\eqref{eq:branch-energy-asymptotic}, we obtain
\begin{equation}  \label{eq:comp_W} 
     W(T)
    \le
    rn\sum_{i=1}^k t_i
    -\sum_{i=1}^k\bigl(b_i n+R_i t_i^2\bigr)
    +O_r(n).
\end{equation}
For every $i$, by AM-GM,
\[
    b_i n+R_i t_i^2
    \ge 2t_i\sqrt{b_iR_i n}
    \ge 2t_i\sqrt{\mu_r n}.
\]
Since $\sum_i t_i=n-1$, it follows that
\begin{align*}
    W(T)
    &\le
    rn(n-1)-2\sqrt{\mu_r n}\,(n-1)+O_r(n)\\
    &=
    rn^2-2\sqrt{\mu_r}\,n^{3/2}+O_r(n).
\end{align*}

For the reverse inequality, choose a skeleton $S$ attaining $\mu_r$ and
write
\[
    b:=b(S),
    \qquad
    R:=R(S),
    \qquad
    bR=\mu_r.
\]
Let
\[
    k=\sqrt{\frac{Rn}{b}}+O(1).
\]
Divide $n-1$ into $k$ positive integers $t_1,\ldots,t_k$ differing by at most
one. Thus
\[
    \sum_{i=1}^k t_i^2=\frac{(n-1)^2}{k}+O(k).
\]
Construct $k$ branches with skeleton $S$ and orders $t_i$, distributing
their pendant leaves according to a minimizing vector in the definition of
$R(S)$. 

Choose the pendant-leaf multiplicities of each branch so that the upper estimate for $\Phi(B_i)$ in~\eqref{eq:branch-energy-asymptotic} is attained up to $O_r(t_i)$. By the exact edge formula and~\eqref{eq:branch-linear-term}, and using the uniformity of the error term in~\eqref{eq:branch-energy-asymptotic}, we obtain
	\begin{align*}
		W(T)&=rn(n-1)-kbn-R\sum_{i=1}^k t_i^2+O_r(n)\\
		&=rn^2-kbn-\frac{Rn^2}{k}+O_r(n)\\
		&=rn^2-2\sqrt{\mu_r}\,n^{3/2}+O_r(n).
	\end{align*}
	Here we used $\sum_i t_i=n-1$, $\sum_i t_i^2=(n-1)^2/k+O(k)$ and $k=\sqrt{Rn/b}+O(1)$. In particular, the accumulated rounding error remains $O_r(n)$.

When $\diam(T)=2r+1$, let $c_1c_2$ be the central edge, and let $x$ be the order of the component containing $c_1$ after deleting that edge. The total order of all branches attached to $c_1$ and $c_2$ is $N=n-2$. The central edge contributes
	\[
	x(n-x)=\frac{n^2}{4}-\left(x-\frac n2\right)^2.
	\]
	The same edge formula and branch estimates as above give
	\[
	W(T)\leq rnN+\frac{n^2}{4}-2N\sqrt{\mu_r n}+O_r(n)
	=\left(r+\frac14\right)n^2-2\sqrt{\mu_r}\,n^{3/2}+O_r(n).
	\]
	For the reverse inequality, choose a skeleton $S$ attaining $\mu_r=bR$. Split $N$ into two integers $N_1,N_2$ differing by at most one. For $j=1,2$, choose $k_j=N_j/\sqrt{bn/R}+O(1)$ branches of skeleton $S$ attached to $c_j$, with orders differing by at most one and summing to $N_j$. Every branch order is $\sqrt{bn/R}+O_r(1)$. Choose pendant-leaf multiplicities as in (3). If $k=k_1+k_2$, then $k=\sqrt{Rn/b}+O_r(1)$ and $\sum_i t_i^2= N^2/k+O_r(k)$, since all branch orders differ from their common mean by $O_r(1)$. The exact edge formula yields
	\[
	W(T)=rnN+\frac{n^2}{4}-kbn-R\frac{N^2}{k}+O_r(n)
	=\left(r+\frac14\right)n^2-2\sqrt{\mu_r}\,n^{3/2}+O_r(n).
	\]
	The central-edge term differs from $n^2/4$ by at most $1/4$ in this construction. This completes the proof of the odd-diameter case.
    
\end{proof}

The numbers $\mu_r$ can be computed by a finite bottom-up recursion. Let
\[
    \mathcal P_0=\{(1,1)\}.
\]
For $j\ge 1$, form
\begin{align*}
    \mathcal Q_j={}&
    \bigl\{(b+j+1,R+1):(b,R)\in\mathcal P_{j-1}\bigr\}\\
    &\mathbin{\cup}
    \left\{
    \left(
        b_1+b_2+j+1,\,
        1+\frac{R_1R_2}{R_1+R_2}
    \right):
    (b_i,R_i)\in\mathcal P_{j-1}
    \right\}.
\end{align*}

The first operation corresponds to adding a root with one child. The second
corresponds to adding a root with two children: the two child resistances are
combined (with the optimal ratio $x$ of the number of leaves in the first part with respect to the total) in parallel and then placed in series with the new attachment
edge.
Here the only non-trivial computation is that $x^2R_1+(1-x)^2R_2= x^2(R_1+R_2)-2R_2x+R_2$ is minimized when $x= \frac{R_2}{R_1+R_2}$, in which case $x^2R_1+(1-x)^2R_2=\frac{R_1R_2}{R_1+R_2}.$ 

Let $\mathcal P_j$ be obtained from $\mathcal Q_j$ by deleting every
coordinatewise dominated pair: $(b,R)$ is deleted if there is another pair
$(b',R')$ satisfying
\[
    b'\le b,\qquad R'\le R,
\]
with at least one strict inequality.

Since both operations are increasing in both coordinates, deleting
dominated pairs does not affect any later minimum. Hence
\[
    \mu_r=\min\{bR:(b,R)\in\mathcal P_{r-2}\}.
\]

Applying this recursion with exact rational arithmetic gives the following
values\footnote{The corresponding program is also made publicly available at\\ \url{https://github.com/JorikJooken/maximizeTotalDistanceGivenOrderAndRadius}.}. The middle two columns display one minimizing pair.
\[
\begin{array}{c|c|c|c}
 r & b(S) & R(S) & \mu_r \\ \hline
 3  & 3  & 2       & 6 \\
 4  & 7  & 5/2     & 35/2 \\
 5  & 11 & 7/2     & 77/2 \\
 6  & 16 & 9/2     & 72 \\
 7  & 24 & 5       & 120 \\
 8  & 31 & 6       & 186 \\
 9  & 39 & 7       & 273 \\
 10 & 51 & 15/2    & 765/2 \\
 11 & 61 & 17/2    & 1037/2 \\
 12 & 78 & 35/4    & 1365/2 \\
 13 & 90 & 39/4    & 1755/2
\end{array}
\]

\begin{figure}
    \centering
    \begin{tikzpicture}[scale=1,
every node/.style={circle,fill,inner sep=1.5pt}]

\node (a0) at (0,0) {};
\node (a1) at (0,-1) {};
\node (a2) at (0,-2) {};
\node (a3) at (0,-3) {};

\draw (a0)--(a1)--(a2)--(a3);

\end{tikzpicture}
\quad
\begin{tikzpicture}[scale=1,
every node/.style={circle,fill,inner sep=1.5pt}]

\node (a1) at (0,-1) {};
\node (a2) at (0,-2) {};
\node (b) at (-0,-3) {};
\draw (a1)--(a2)--(b);

\node (c1) at (-1,-4) {};
\node (c2) at ( 1,-4) {};

\draw (b)--(c1);
\draw (b)--(c2);

\end{tikzpicture}
\quad
\begin{tikzpicture}[scale=1,
every node/.style={circle,fill,inner sep=1.5pt}]

\node (a1) at (0,-1) {};
\node (a2) at (0,-2) {};

\draw (a1)--(a2);

\node (b1) at (-0.5,-3) {};
\node (b2) at ( 0.5,-3) {};

\draw (a2)--(b1);
\draw (a2)--(b2);

\node (c1) at (-1,-4) {};
\node (c2) at ( 1,-4) {};

\draw (b1)--(c1);
\draw (b2)--(c2);

\end{tikzpicture}\quad
\begin{tikzpicture}[scale=1,
every node/.style={circle,fill,inner sep=1.5pt}]

\node (a2) at (0,-2) {};
\node (a3) at (0,-3) {};

\draw (a2)--(a3);

\node (b1) at (-1,-4) {};
\node (b2) at ( 1,-4) {};

\draw (a3)--(b1);
\draw (a3)--(b2);

\node (c1) at (-1.5,-5) {};
\node (c2) at (-0.5,-5) {};
\node (c3) at ( 1,-5) {};

\draw (b1)--(c1);
\draw (b1)--(c2);
\draw (b2)--(c3);

\end{tikzpicture}
\quad
\quad
\begin{tikzpicture}[scale=1,
every node/.style={circle,fill,inner sep=1.5pt}]

\node (a2) at (0,-2) {};
\node (a3) at (0,-3) {};

\draw (a2)--(a3);

\node (b1) at (-1,-4) {};
\node (b2) at ( 1,-4) {};

\draw (a3)--(b1);
\draw (a3)--(b2);

\node (c1) at (-1.5,-5) {};
\node (c2) at (-0.5,-5) {};
\node (c3) at ( 0.5,-5) {};
\node (c4) at ( 1.5,-5) {};

\draw (b1)--(c1);
\draw (b1)--(c2);

\draw (b2)--(c3);
\draw (b2)--(c4);





\end{tikzpicture}

    \caption{5 possible skeletons $\widehat S$ for the branches of a tree with radius $4$}
    \label{fig:5Skeletons_r4}
\end{figure}
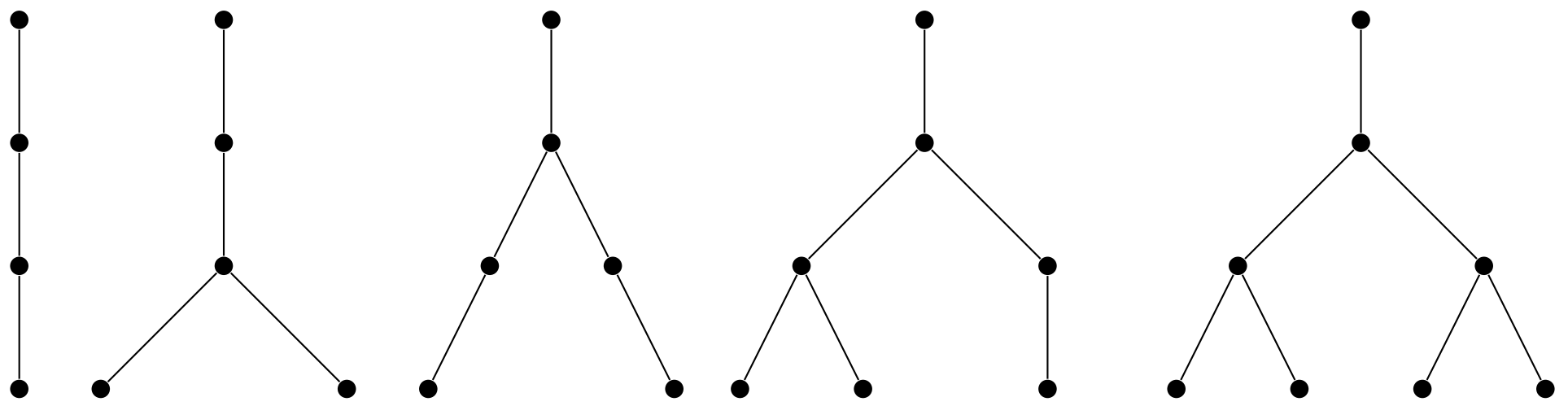

\begin{remark}
For $3\le r\le 11$, a minimizing skeleton can be chosen to have only one
branching vertex. This pattern first fails at $r=12$. 

More precise computations verify that every branch of the extremal tree for $r=3$ is a broom (so only one of the two skeletons appear), while for $r=4$ none of them is when $n>597$. In the latter case, all skeletons of the branches are of the second form as presented in~\cref{fig:5Skeletons_r4}, while for some $n \le 597$ multiple types of skeletons appear in the extremal tree.
For general fixed $r$, \cref{sec:app_mu_r_asymptotic_sharpness}
shows that, for all sufficiently large $n$, every branch of every
tree of maximum total distance of order $n$ and diameter $2r$ or $2r+1$
has a skeleton attaining $\mu_r$.

\end{remark}

\section{Structure of the extremal trees}\label{sec:structure}

Given~\cref{sec:given_r}, it is natural to expect that the extremal tree is very balanced. In this section, we want to stress that this is not true for large $r$ and moderate values of $n$. The reason for this is that various skeletons may appear in an extremal tree, and as such lead to less balanced branches. On the other hand, we also prove that branches with the same skeleton and attached to the same central vertex indeed have orders that differ by no more than one.

\begin{figure}
\begin{tikzpicture}[
    every node/.style={circle,fill,inner sep=1.5pt}
]


\node (a1) at (0,-1) {};
\node (a2) at (0,-2) {};
\node (a3) at (0,-3) {};

\draw (a1)--(a2)--(a3);

\node[draw=none,fill=none] at (0,-3.375) {$\vdots$};

\node (a9)  at (0,-4) {};
\node ()  at (0,-5) {};
\node (a10) at (0,-6) {};
\node (b)   at (0,-7) {};

\draw (a9)--(a10)--(b);


\node (c1) at (-2.5,-8) {};
\node (c2) at ( 2.5,-8) {};
\node () at (-2.5,-9) {};
\node () at ( 2.5,-9) {};

\draw (b)--(c1);
\draw (b)--(c2);

\node (d1) at (-2.5,-10) {};
\node (d2) at ( 2.5,-10) {};

\draw (c1)--(d1);
\draw (c2)--(d2);


\node (e1) at (-3.5,-11) {};
\node (e2) at (-1.5,-11) {};
\node (e3) at ( 1.5,-11) {};
\node (e4) at ( 3.5,-11) {};

\draw (d1)--(e1);
\draw (d1)--(e2);

\draw (d2)--(e3);
\draw (d2)--(e4);


\node (f1) at (-3.95,-12) {};
\node (f2) at (-3.75,-12) {};
\node (f3) at (-3.25,-12) {};
\node (f4) at (-3.05,-12) {};

\draw (e1)--(f1);
\draw (e1)--(f2);
\draw (e1)--(f3);
\draw (e1)--(f4);

\node[draw=none,fill=none] at (-3.5,-12) {$\cdots$};

\node[draw=none,fill=none] at (-3.5,-12.4) {$47$ leaves};


\node (g1) at (-1.95,-12) {};
\node (g2) at (-1.75,-12) {};
\node (g3) at (-1.25,-12) {};
\node (g4) at (-1.05,-12) {};

\draw (e2)--(g1);
\draw (e2)--(g2);
\draw (e2)--(g3);
\draw (e2)--(g4);

\node[draw=none,fill=none] at (-1.5,-12) {$\cdots$};

\node[draw=none,fill=none] at (-1.5,-12.4) {$47$ leaves};


\node (h1) at (1.05,-12) {};
\node (h2) at (1.25,-12) {};
\node (h3) at (1.75,-12) {};
\node (h4) at (1.95,-12) {};

\draw (e3)--(h1);
\draw (e3)--(h2);
\draw (e3)--(h3);
\draw (e3)--(h4);

\node[draw=none,fill=none] at (1.5,-12) {$\cdots$};

\node[draw=none,fill=none] at (1.5,-12.4) {$47$ leaves};


\node (i1) at (3.05,-12) {};
\node (i2) at (3.25,-12) {};
\node (i3) at (3.75,-12) {};
\node (i4) at (3.95,-12) {};

\draw (e4)--(i1);
\draw (e4)--(i2);
\draw (e4)--(i3);
\draw (e4)--(i4);

\node[draw=none,fill=none] at (3.5,-12) {$\cdots$};

\node[draw=none,fill=none] at (3.5,-12.4) {$\textbf{47/48}$ leaves};



\end{tikzpicture}
    \centering
\begin{tikzpicture}[
    every node/.style={circle,fill,inner sep=1.5pt}, scale=0.9
]


\node (a1) at (0,-1) {};
\node (a2) at (0,-2) {};
\node (a3) at (0,-3) {};

\draw (a1)--(a2)--(a3);
\node[draw=none,fill=none] at (0,-3.375) {$\vdots$};

\node (a9)  at (0,-4) {};
\node ()  at (0,-5) {};
\node (a10) at (0,-6) {};
\node (b)   at (0,-7) {};

\draw (a9)--(a10)--(b);


\node (c1) at (-2.5,-8) {};
\node (c2) at ( 2.5,-8) {};
\node () at (-2.5,-9) {};
\node () at ( 2.5,-9) {};

\draw (b)--(c1);
\draw (b)--(c2);

\node (d1) at (-2.5,-10) {};
\node (d2) at ( 2.5,-10) {};

\draw (c1)--(d1);
\draw (c2)--(d2);


\node (e1) at (-3.5,-11) {};
\node (e2) at (-1.5,-11) {};
\node (e3) at ( 2.5,-11) {};

\draw (d1)--(e1);
\draw (d1)--(e2);

\draw (d2)--(e3);


\node (f1) at (-3.95,-12) {};
\node (f2) at (-3.75,-12) {};
\node (f3) at (-3.25,-12) {};
\node (f4) at (-3.05,-12) {};

\draw (e1)--(f1);
\draw (e1)--(f2);
\draw (e1)--(f3);
\draw (e1)--(f4);

\node[draw=none,fill=none] at (-3.5,-12) {$\cdots$};

\node[draw=none,fill=none] at (-3.5,-12.4) {$50$ leaves};


\node (g1) at (-1.95,-12) {};
\node (g2) at (-1.75,-12) {};
\node (g3) at (-1.25,-12) {};
\node (g4) at (-1.05,-12) {};

\draw (e2)--(g1);
\draw (e2)--(g2);
\draw (e2)--(g3);
\draw (e2)--(g4);

\node[draw=none,fill=none] at (-1.5,-12) {$\cdots$};

\node[draw=none,fill=none] at (-1.5,-12.4) {$50$ leaves};


\node (h1) at (1.55,-12) {};
\node (h2) at (1.75,-12) {};
\node (h3) at (3.25,-12) {};
\node (h4) at (3.45,-12) {};

\draw (e3)--(h1);
\draw (e3)--(h2);
\draw (e3)--(h3);
\draw (e3)--(h4);

\node[draw=none,fill=none] at (2.5,-12) {$\cdots$};

\node[draw=none,fill=none] at (2.5,-12.4) {$88$ leaves};



\end{tikzpicture}
\quad
\begin{tikzpicture}[
    every node/.style={circle,fill,inner sep=1.5pt}, scale=0.9
]


\node (a1) at (0,-1) {};
\node (a2) at (0,-2) {};
\node (a3) at (0,-3) {};

\draw (a1)--(a2)--(a3);

\node[draw=none,fill=none] at (0,-3.375) {$\vdots$};

\node (a9)  at (0,-4) {};
\node (a10) at (0,-5) {};
\node (b)   at (0,-6) {};

\draw (a9)--(a10)--(b);


\node (c1) at (-2.5,-7) {};
\node (c2) at ( 2.5,-7) {};
\node () at (-2.5,-8) {};
\node () at ( 2.5,-8) {};
\node () at (-2.5,-9) {};
\node () at ( 2.5,-9) {};

\draw (b)--(c1);
\draw (b)--(c2);

\node (d1) at (-2.5,-10) {};
\node (d2) at ( 2.5,-10) {};

\draw (c1)--(d1);
\draw (c2)--(d2);


\node (e1) at (-3.5,-11) {};
\node (e2) at (-1.5,-11) {};
\node (e3) at ( 1.5,-11) {};
\node (e4) at ( 3.5,-11) {};

\draw (d1)--(e1);
\draw (d1)--(e2);

\draw (d2)--(e3);
\draw (d2)--(e4);


\node (f1) at (-3.95,-12) {};
\node (f2) at (-3.75,-12) {};
\node (f3) at (-3.25,-12) {};
\node (f4) at (-3.05,-12) {};

\draw (e1)--(f1);
\draw (e1)--(f2);
\draw (e1)--(f3);
\draw (e1)--(f4);

\node[draw=none,fill=none] at (-3.5,-12) {$\cdots$};

\node[draw=none,fill=none] at (-3.5,-12.4) {$49$ leaves};


\node (g1) at (-1.95,-12) {};
\node (g2) at (-1.75,-12) {};
\node (g3) at (-1.25,-12) {};
\node (g4) at (-1.05,-12) {};

\draw (e2)--(g1);
\draw (e2)--(g2);
\draw (e2)--(g3);
\draw (e2)--(g4);

\node[draw=none,fill=none] at (-1.5,-12) {$\cdots$};

\node[draw=none,fill=none] at (-1.5,-12.4) {$49$ leaves};


\node (h1) at (1.05,-12) {};
\node (h2) at (1.25,-12) {};
\node (h3) at (1.75,-12) {};
\node (h4) at (1.95,-12) {};

\draw (e3)--(h1);
\draw (e3)--(h2);
\draw (e3)--(h3);
\draw (e3)--(h4);

\node[draw=none,fill=none] at (1.5,-12) {$\cdots$};

\node[draw=none,fill=none] at (1.5,-12.4) {$49$ leaves};


\node (i1) at (3.05,-12) {};
\node (i2) at (3.25,-12) {};
\node (i3) at (3.75,-12) {};
\node (i4) at (3.95,-12) {};

\draw (e4)--(i1);
\draw (e4)--(i2);
\draw (e4)--(i3);
\draw (e4)--(i4);

\node[draw=none,fill=none] at (3.5,-12) {$\cdots$};

\node[draw=none,fill=none] at (3.5,-12.4) {$49$ leaves};



\end{tikzpicture}
 \caption{Branches of order $209/210$, $208$ and $218$ in extremal tree for $r=16, n=4300$}
    \label{fig:subtrees_n4300_r16}
\end{figure}
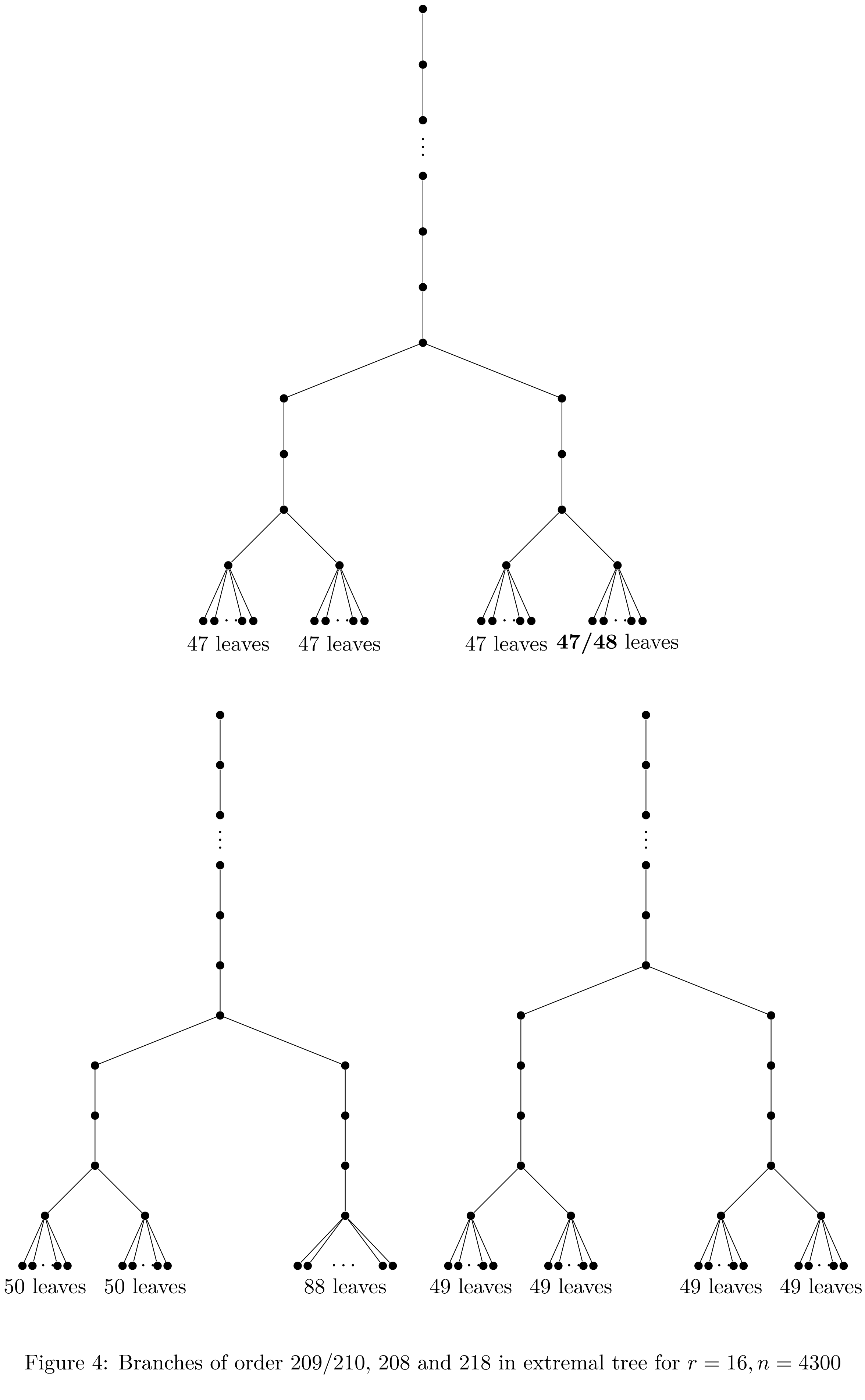

As a concrete example to indicate how hard it can be to predict an extremal tree maximizing $W$ given order $n$ and diameter $d$ or radius $r$, for $n=4300$ and $d=32, r=16$, we present an extremal tree whose branches have $3$ different skeletons. It has a central vertex connected with $2$ branches of order $208$, $1$ branch of order $209$, $4$ branches of order $210$ and $13$ branches of order $218$. Here the branches are depicted in~\cref{fig:subtrees_n4300_r16}.
We notice that the branching does not always happen in the same layer (here both $L_{5}$ and $L_{6}$), the number of branching points in different branch skeletons can be different (here $2$ or $3$) and different branches can have very different order (notice e.g. the $88+4=92$ compared with $2 \cdot 50+5=105$ within the branch of order $208$) and there can be multiple different orders (e.g., here $4$, with $208,209,210$ and $218$).

For odd diameter, it can be observed that often the central edge does not split the tree exactly in half. E.g., for $d=7$ and $n=248$, the two parts have order $121$ and $127$.

\begin{figure}[h]
    \centering
    \begin{tikzpicture}[
    every node/.style={circle,fill,inner sep=1.5pt}, scale=0.89
]


\node (a1) at (0,-1) {};
\node (a2) at (0,-2) {};
\node (a3) at (0,-3) {};

\draw (a1)--(a2)--(a3);

\node[draw=none,fill=none] at (0,-3.375) {$\vdots$};

\node (a9)  at (0,-4) {};
\node (a10) at (0,-5) {};
\node (b)   at (0,-6) {};

\draw (a9)--(a10)--(b);


\node (c1) at (-2.5,-7) {};
\node (c2) at ( 2.5,-7) {};
\node (d1) at (-2.5,-8) {};
\node (d2) at ( 2.5,-8) {};
\node (e1) at (-2.5,-9) {};
\node (e2) at ( 2.5,-9) {};

\node[draw=none,fill=none] at (-2.5,-8.5) {$\vdots$};
\node[draw=none,fill=none] at (2.5,-8.5) {$\vdots$};

\draw (b)--(c1);
\draw (b)--(c2);

\node (f1) at (-2.5,-10) {};
\node (f2) at ( 2.5,-10) {};
\node (g) at ( 2.5,-11) {};

\draw (c1)--(d1);
\draw (c2)--(d2);
\draw (e1)--(f1);
\draw (e2)--(f2);


\node (g1) at (-3.5,-11) {};
\node (g2) at (-1.5,-11) {};
\node (h3) at ( 1.5,-12) {};
\node (h4) at ( 3.5,-12) {};

\draw (f1)--(g1);
\draw (f1)--(g2);
\draw (f2)--(g);
\draw (g)--(h3);
\draw (g)--(h4);

\node () at (-3.5,-12) {};
\node () at (-1.5,-12) {};

\node () at (-3.5,-13) {};
\node () at (-1.5,-13) {};

\node () at (-3.5,-14) {};
\node (r1) at (-2,-14) {};
\node (r2) at (-1,-14) {};

\draw (g1)--(-3.5,-14) ;
\draw (g2)-- (-1.5,-13)--(r1);
\draw  (-1.5,-13)--(r2);

\node (i1) at ( 1.5,-13) {};
\node () at ( 3.5,-13) {};
\node (j1) at ( 1,-14) {};
\node (j2) at ( 2,-14) {};
\node (i2) at ( 3.5,-14) {};

\draw (h4)--(i2);
\draw (i1)--(h3);
\draw (j1)--(i1)--(j2);

\end{tikzpicture}
 \caption{Optimal skeleton (for $n$ large) for $r=45$: splits are at level $43,41,40$ and $28$}
    \label{fig:skeleton_r45}
\end{figure}
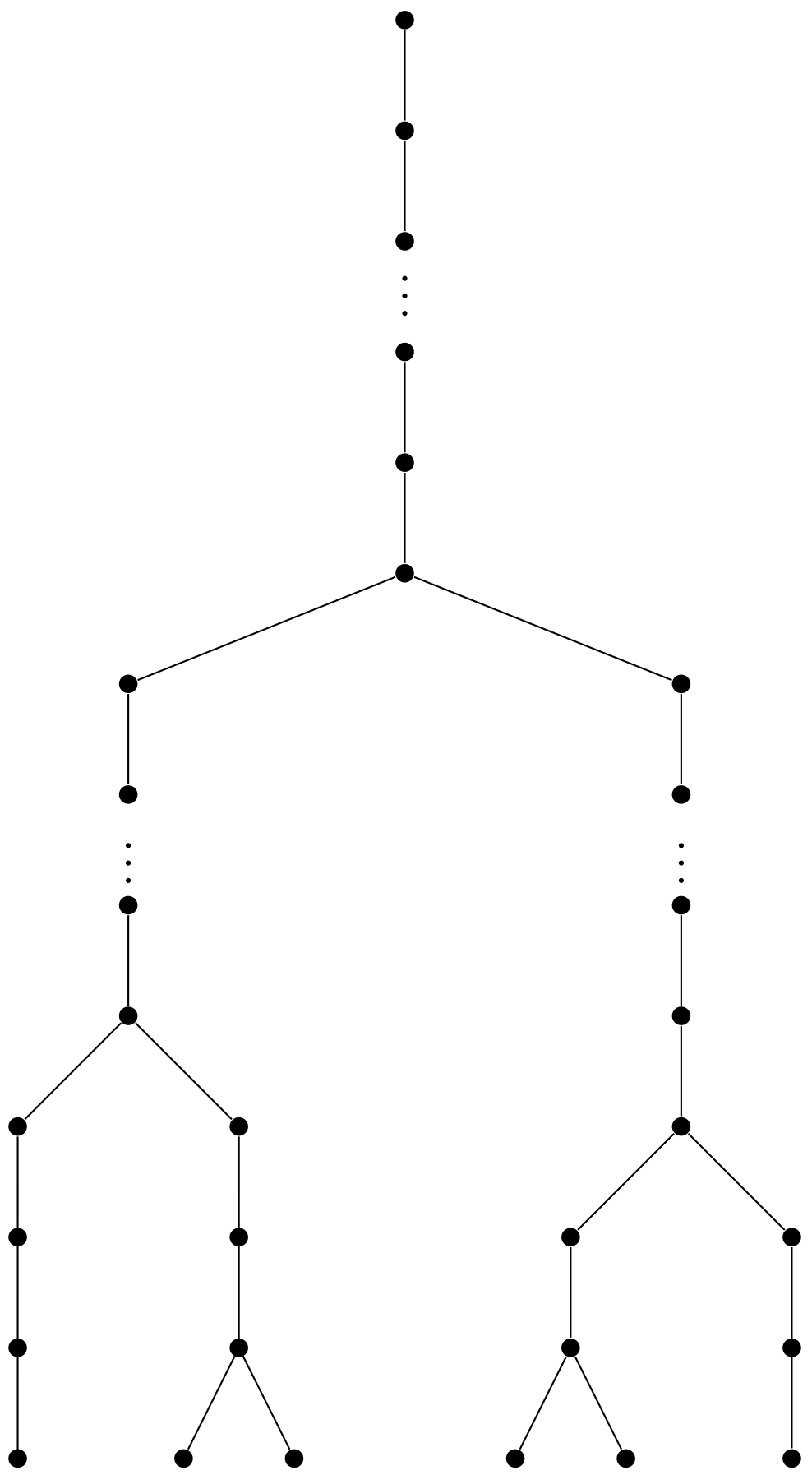

It can be expected that the skeleton(s) attaining $\mu_r$ are always symmetric.
This turns out to be not always the case. 
E.g. for $r=45$, the optimal skeleton is presented in~\cref{fig:skeleton_r45}.

Finally, we prove the balancedness among branches with the same skeleton at the same side.

	\begin{prop}
	   Let $T$ be a tree maximizing the Wiener index among trees of its order and diameter. If two branches attached to the same central vertex have isomorphic skeletons, then their orders differ by at most one.  
	\end{prop}
    
	\begin{proof}
		Let $B_1$ and $B_2$ be two such branches attached to a central vertex $c$, and identify their skeletons with a common rooted skeleton $S$. Enumerate the terminal vertices $z_1,\ldots,z_q$ of $S$ in depth-first order. Let $a_j$ and $b_j$ be the numbers of pendant leaves attached to $z_j$ in the two copies of $S$, respectively. Suppose, without loss of generality, that
		\[
		|V(B_1)|-|V(B_2)|=\sum_{j=1}^q(a_j-b_j)\geq2.
		\]
		We construct two new branches with skeleton $S$ by redistributing the pendant leaves. Form a sequence consisting of $a_1+b_1$ tokens of type $1$, followed by $a_2+b_2$ tokens of type $2$, and so on. Allocate the tokens alternately to the first and second new branches. For each $j$, let $a'_j$ and $b'_j$ be the numbers of tokens of type $j$ allocated to the respective branches, and attach the corresponding numbers of pendant leaves to their copies of $z_j$. We have $a'_j+b'_j=a_j+b_j$, and since $a_j+b_j\geq2$, both $a'_j$ and $b'_j$ are positive. Thus the skeletons and branch heights are preserved, as are the order and diameter of the whole tree.
		
		For any edge $e$ of $\widehat S$, including the attachment edge at $c$, the terminals below $e$ are consecutive in depth-first order. Consequently, the tokens associated with these terminals form a consecutive interval of the full sequence. The numbers allocated to the two new branches in this interval differ by at most one. Let $s_1(e)$ and $s_2(e)$ denote the orders of the components below the two original copies of $e$, and let $s'_1(e)$ and $s'_2(e)$ be the analogous orders after redistribution. The numbers of skeleton vertices below corresponding edges are equal, so
		\[
		s'_1(e)+s'_2(e)=s_1(e)+s_2(e),\qquad
		|s'_1(e)-s'_2(e)|\leq1.
		\]
		Put $f(x)=x(n-x)$, where $n=|V(T)|$. Since $f$ is strictly concave, for integers with a fixed sum the quantity $f(x)+f(y)$ is maximized when $|x-y|\leq1$. In fact, for a fixed sum $m=x+y$,
		\[
		f(x)+f(y)=nm-\frac{m^2+(x-y)^2}{2}.
		\]
		It follows that the combined contribution of every pair of corresponding skeleton edges does not decrease. For the attachment edges at $c$, the original component orders differ by at least two, so this increase is strict.
		
		The combined contribution of pendant edges is unchanged, because their total number is unchanged and every pendant edge contributes $n-1$. For any edge outside the two branches, the two branches lie on the same side of the edge, since they are attached to the same vertex $c$. Their combined order is unchanged, and hence the contribution of every such edge is unchanged. By the edge formula for the Wiener index, the new tree therefore has strictly larger Wiener index than $T$, a contradiction.
	\end{proof}

We observe that an extremal graph has at most one block not isomorphic to $K_2$.
This is consistent with~\cref{conj:Nlarge}.

\begin{prop}
    A graph maximizing the total distance among connected graphs of order $n$ and diameter $d$ has at most one block not isomorphic to $K_2$. Equivalently, it is either a tree or is obtained from a $2$-connected graph by attaching trees at its vertices.
\end{prop}

\begin{proof}
Suppose for the sake of obtaining a contradiction that an extremal graph $G$ has two blocks not isomorphic to $K_2$. Let $v$ be a cutvertex and let $A$ and $B$ be two components of $G-v$ such that both $G[A\cup\{v\}]$ and $G[B\cup\{v\}]$ contain a cycle. Let $a,b$ be the maximum distances from $v$ in $G[A\cup\{v\}]$ and $G[B\cup\{v\}]$, respectively. Since $a+b\le d$, we may assume that $2a\le d$.

Replace $G[A\cup\{v\}]$ by a breadth-first spanning tree rooted at $v$ (i.e., delete all edges in $G[A\cup\{v\}]$ that do not belong to this breadth-first spanning tree). This does not modify any distance from $v$ to a vertex in $G[A\cup\{v\}]$ and therefore all distances involving vertices outside $A$ are unchanged. Distances within the replaced subgraph are at most $2a\le d$ and therefore the resulting spanning subgraph still has diameter $d$. Since at least one edge is deleted, some distances increase strictly and none decrease, contradicting the maximality of $W(G)$.
\end{proof}

\section{Conclusion}\label{sec:conc}


We described an algorithm that determines the maximum total distance
and constructs an extremal graph for given order $n$ and
radius $r$ in $O(n^2r)$ time, thereby resolving the algorithmic
version of the problem. We complement this algorithm with a precise asymptotic estimate for the maximum total distance when the radius is fixed, and
with results on the structure of extremal trees.

Our examples show that extremal trees can exhibit substantial
structural variation: branches of a single extremal tree may
have different skeletons and different orders. This variability
suggests that a simple closed-form characterization valid for
all $n$ and $r$ should not be expected. An exact algorithm,
complemented by asymptotic and structural results, therefore
provides a natural form of resolution of the problem.

When studying the same problem, but for given diameter $d$, as posed by Plesn\'{\i}k, we now have a reasonably good understanding of the extremal behaviour. In~\cite{SIRV19}, it was observed for $d=4$ that the extremal graphs for the problem of Plesn\'{\i}k should be trees except when $n=9$. It is natural to think this may be true for larger even diameters as well, but the intermediate range of exceptions is easily noted to be larger.
For example, for $n=2d$ we can compute the extremal tree and compare with $C_{2d}$ and note that for large $d$, the maximum Wiener index given order $n=2d$ and diameter $d$ is not obtained by a tree.

The results suggest a qualitative picture in which different
graph structures are extremal in different ranges of the order
relative to the diameter. A possibility is as follows. For $n$ small, the extremal graph is a tree. This is easily seen to be true for $n<3/2d$, as any such graph has a diameter-preserving spanning tree, in which case our algorithms can compute the optimum in $O(d^3)$.
For $n$ around $2d$, it is approximately a cycle.
For $n$ large, it is expected to be as mentioned in~\cref{conj:Nlarge}, with the structure depending on the parity of $d$.

When the extremal graph is a tree or clique-rooted tree, we have constructed algorithms in~\cref{sec:algo} that can compute the exact maximum and the extremal graphs of the latter form.
As such, the main (challenging) problem that remains for the future, is exactly~\cref{conj:Nlarge}, which is concerned with the extremal structure for large order.

\section*{AI Declaration}

At the end of this project, ChatGPT 5.6 Sol has been used for further discussion and checking. This led to the idea for~\cref{thr:fixed-diameter-asymptotics}, with the precise version being put in appendix. The main manuscript is the work of the authors.


\appendix

\crefalias{section}{appendix}
\section{Stronger form of~\cref{thr:fixed-diameter-asymptotics}}\label{sec:app_mu_r_asymptotic_sharpness}

\begin{lem}[Second-order expansion of the branch energy]
\label{lem:second-order-branch-energy}
Let $S$ be a fixed $r$-branch skeleton. For every feasible integer $t$, let
\[
    H_S(t):=
    \min_B\sum_{e\in E^+(B)}s_e^2,
\]
where the minimum is over all branches $B$ of order $t$ with skeleton $S$,
the set $E^+(B)$ consists of the edges of $B$ together with its attachment
edge, and $s_e$ is the order of the component below $e$. Then there exists
a constant $\gamma(S)$ such that
\[
    H_S(t)=R(S)t^2+\gamma(S)t+O_S(1).
\]
In particular, since the collection of $r$-branch skeletons is finite, the
error term may be taken to be $O_r(1)$ uniformly over all such skeletons.
\end{lem}

\begin{proof}
Let $Z(S)=\{z_1,\ldots,z_q\}$ be the set of terminal vertices of $S$, and
let $a_j\ge1$ be the number of pendant leaves attached to $z_j$. Put
\[
    A:=\sum_{j=1}^q a_j=t-|V(S)|.
\]
For each edge $e$ of $\widehat S$, let $p_e$ be the number of skeleton
vertices below $e$, and let
\[
    M_{ej}=
    \begin{cases}
        1,&\text{if $z_j$ lies below $e$},\\
        0,&\text{otherwise}.
    \end{cases}
\]
Writing $\mathbf a=(a_1,\ldots,a_q)^{\mathsf T}$ and
$\mathbf p=(p_e)_e$, the quadratic contribution of the branch is
\[
    \Phi_S(\mathbf a)
    =
    A+\|M\mathbf a+\mathbf p\|_2^2.
\]
The term $A$ comes from the $A$ pendant edges, each of which has
$s_e=1$.

First minimize $\Phi_S$ over real vectors satisfying
\[
    \mathbf 1^{\mathsf T}\mathbf a=A.
\]
Since the quadratic form $\|M\mathbf a\|_2^2$ is positive definite, the
minimizer is unique. The Lagrange multiplier equations show that it has
the form
\[
    \mathbf a^*(A)=A\boldsymbol\lambda+\boldsymbol\xi,
\]
where $\boldsymbol\lambda$ and $\boldsymbol\xi$ depend only on $S$.
Moreover,
\[
    \mathbf 1^{\mathsf T}\boldsymbol\lambda=1
\]
and $\boldsymbol\lambda$ minimizes
\[
    \|M\boldsymbol\lambda\|_2^2
\]
subject to this constraint. Hence
\[
    \|M\boldsymbol\lambda\|_2^2=R(S).
\]
Every coordinate of $\boldsymbol\lambda$ is positive: in the electrical
interpretation, a positive amount of the unit current reaches every
terminal of the finite tree. Thus the positivity constraints
$a_j\ge1$ are inactive for all sufficiently large $A$.

The minimum over the real affine hyperplane is therefore a quadratic
polynomial
\[
    R(S)A^2+\ell(S)A+c(S).
\]
For every integer $A$, one can choose an integral vector
$\widetilde{\mathbf a}$ with
\[
    \mathbf 1^{\mathsf T}\widetilde{\mathbf a}=A,
    \qquad
    \|\widetilde{\mathbf a}-\mathbf a^*(A)\|_\infty=O_S(1).
\]
Since $\mathbf a^*(A)$ is the constrained real minimizer, the linear term
in the difference vanishes along the hyperplane
$\mathbf 1^{\mathsf T}\mathbf a=A$. Consequently,
\[
    \Phi_S(\widetilde{\mathbf a})
    -
    \Phi_S(\mathbf a^*(A))
    =
    O_S(1).
\]
The integral minimum is bounded below by the real minimum, so
\[
    H_S(t)=R(S)A^2+\ell(S)A+c(S)+O_S(1).
\]
Substituting $A=t-|V(S)|$ gives
\[
    H_S(t)=R(S)t^2+\gamma(S)t+O_S(1)
\]
for a suitable constant $\gamma(S)$. Enlarging the bounded error accounts
for the finitely many small feasible values of $t$.
\end{proof}

\begin{prop}[Only minimizing skeletons occur]
\label{prop:only-minimizing-skeletons}
Fix $r\ge2$, and let
\[
    \mathcal M_r:=
    \bigl\{
        S:\ b(S)R(S)=\mu_r
    \bigr\}.
\]
There exists $n_0(r)$ such that, for every $n\ge n_0(r)$, every tree
maximizing the Wiener index among the trees of order $n$ and diameter
$2r$, or among the trees of order $n$ and diameter $2r+1$, has the
following property:
\[
    S\in\mathcal M_r
\]
for the skeleton $S$ of every branch.

In particular, if the skeleton attaining $\mu_r$ is unique, then every
branch of every sufficiently large extremal tree has that skeleton.
\end{prop}

\begin{proof}
For a skeleton $S$, abbreviate
\[
    b_S:=b(S),\qquad R_S:=R(S),\qquad \gamma_S:=\gamma(S).
\]
Choose
\[
    \gamma_0:=
    \min\{\gamma_S:S\in\mathcal M_r\}
\]
and define
\[
    \Lambda_n:=2\sqrt{\mu_r n}+\gamma_0.
\]
For a branch of order $t$ and skeleton $S$, define its reduced cost by
\[
    G_{S,n}(t)
    :=
    b_Sn+H_S(t)-\Lambda_nt.
\]

Suppose first that $S\in\mathcal M_r$. Since $b_SR_S=\mu_r$,
Lemma~\ref{lem:second-order-branch-energy} gives
\begin{align}
    G_{S,n}(t)
    &=
    b_Sn+R_St^2
    -2\sqrt{\mu_rn}\,t
    +(\gamma_S-\gamma_0)t+O_r(1)\notag\\
    &=
    R_S\left(
        t-\frac{\sqrt{\mu_rn}}{R_S}
    \right)^2
    +(\gamma_S-\gamma_0)t+O_r(1).
    \label{eq:reduced-good}
\end{align}
In particular,
\[
    G_{S,n}(t)
    \ge
    R_S\left(
        t-\frac{\sqrt{\mu_rn}}{R_S}
    \right)^2-C_r.
\]
Thus a branch with $S\in\mathcal M_r$ can have negative reduced cost only
if
\[
    t=\frac{\sqrt{\mu_rn}}{R_S}+O_r(1).
\]
Every such branch has order at least $c_r\sqrt n$, for some $c_r>0$.
Since the sum of all branch orders is at most $n$, there are only
$O_r(\sqrt n)$ branches of negative reduced cost. It follows that the
total reduced cost of all branches whose skeletons lie in
$\mathcal M_r$ is bounded below by
\begin{equation}
    -O_r(\sqrt n).
    \label{eq:total-good-lower}
\end{equation}

Now suppose that $S\notin\mathcal M_r$. Put
\[
    d_S:=\gamma_S-\gamma_0.
\]
Completing the square gives
\begin{align*}
    G_{S,n}(t)
    &=
    R_St^2+
    \bigl(d_S-2\sqrt{\mu_rn}\bigr)t+b_Sn+O_r(1)\\
    &=
    R_S\left(
        t-\frac{2\sqrt{\mu_rn}-d_S}{2R_S}
    \right)^2\\
    &\quad
    +n\left(b_S-\frac{\mu_r}{R_S}\right)
    +\frac{d_S\sqrt{\mu_rn}}{R_S}
    -\frac{d_S^2}{4R_S}
    +O_r(1).
\end{align*}
Since
\[
    b_S-\frac{\mu_r}{R_S}
    =
    \frac{b_SR_S-\mu_r}{R_S}>0,
\]
and there are only finitely many $r$-branch skeletons, there exists
$\varepsilon_r>0$ such that, for all sufficiently large $n$,
\begin{equation}
    G_{S,n}(t)\ge\varepsilon_r n
    \label{eq:bad-linear-loss}
\end{equation}
for every feasible $t$ and every $S\notin\mathcal M_r$.

We next construct a comparison tree whose total reduced cost is only
$O_r(\sqrt n)$. Choose $S_0\in\mathcal M_r$ such that
$\gamma_{S_0}=\gamma_0$, and put
\[
    \tau_n:=\frac{\sqrt{\mu_rn}}{R_{S_0}}.
\]

Consider first diameter $2r$, and put $N=n-1$. Choose an integer
\[
    k=\frac{N}{\tau_n}+O(1)
\]
and divide $N$ into $k$ branch orders $t_1,\ldots,t_k$ differing by at
most one. Then
\[
    t_i=\tau_n+O_r(1)
\]
for every $i$. By \eqref{eq:reduced-good},
\[
    G_{S_0,n}(t_i)=O_r(1).
\]
As $k=O_r(\sqrt n)$, the resulting tree satisfies
\begin{equation}
    \sum_{i=1}^kG_{S_0,n}(t_i)=O_r(\sqrt n).
    \label{eq:comparison-reduced}
\end{equation}
The pendant leaves in every branch are chosen to attain $H_{S_0}(t_i)$.

For an arbitrary tree $T$ of diameter $2r$, let its branch orders and
skeletons be $t_i$ and $S_i$. By the edge formula,
\[
    rn(n-1)-W(T)
    =
    \sum_i\bigl(b_{S_i}n+\Phi(B_i)\bigr),
\]
where $\Phi(B_i)\ge H_{S_i}(t_i)$. Since $\sum_i t_i=n-1$, we obtain
\[
    rn(n-1)-W(T)
    \ge
    \Lambda_n(n-1)+\sum_iG_{S_i,n}(t_i).
\]
If at least one branch has a skeleton outside $\mathcal M_r$, then
\eqref{eq:total-good-lower} and \eqref{eq:bad-linear-loss} give
\[
    rn(n-1)-W(T)
    \ge
    \Lambda_n(n-1)+\varepsilon_rn-O_r(\sqrt n).
\]
On the other hand, the comparison tree from
\eqref{eq:comparison-reduced} satisfies
\[
    rn(n-1)-W(T_0)
    \le
    \Lambda_n(n-1)+O_r(\sqrt n).
\]
For sufficiently large $n$, these inequalities contradict the extremality
of $T$.

For diameter $2r+1$, let $c_1c_2$ be the central edge and let $x$ be the
order of the component containing $c_1$ after that edge is removed. If
$N=n-2$ is the total order of all branches, then
\[
    W(T)
    =
    rnN+x(n-x)
    -
    \sum_i\bigl(b_{S_i}n+\Phi(B_i)\bigr).
\]
Equivalently,
\[
    rnN+\frac{n^2}{4}-W(T)
    =
    \left(x-\frac n2\right)^2
    +
    \sum_i\bigl(b_{S_i}n+\Phi(B_i)\bigr).
\]
The additional square is nonnegative, so the preceding lower bound remains
valid.

For the comparison construction, split $N$ into integers $N_1,N_2$
differing by at most one. Independently partition each $N_j$ into branch
orders equal to $\tau_n+O_r(1)$ and use skeleton $S_0$ for every branch.
The two components of the central edge then differ in order by at most one,
so
\[
    \left(x-\frac n2\right)^2=O(1),
\]
while the total reduced cost is again $O_r(\sqrt n)$. Hence one
non-minimizing skeleton again incurs an $\varepsilon_rn$ loss that cannot
be compensated by the $O_r(\sqrt n)$ rounding error. This proves the
claim.
\end{proof}

\end{document}